\documentclass[12pt,twoside]{amsart}

\author{Florin Ambro}
\address{Department of Mathematical Sciences\\
University of Tokyo,
Komaba, Meguro-Ku,Tokyo 153-8914, JAPAN.}
\email{ambro@ms.u-tokyo.ac.jp}

\newcommand{\1}{{\bf 1}}
\newcommand{\bC}{{\mathbb C}}

\newcommand{\R}{{\mathbb R}}

\newcommand{\cZ}{{\mathcal Z}}

\newcommand{\Exc}{\operatorname{Exc}}

\newcommand{\mult}{\operatorname{mult} }

\newcommand{\relint}{\operatorname{relint}}
\newcommand{\Shed}{\operatorname{Shed}}

\newcommand{\Supp}{\operatorname{Supp}}
\newcommand{\tr}{\operatorname{tr} }

\theoremstyle{plain}
\newtheorem{thm}{Theorem}[section]

\newtheorem{lem}[thm]{Lemma}

\newtheorem{prop}[thm]{Proposition}

\theoremstyle{definition}
\newtheorem{defn}[thm]{Definition}

\newtheorem{ack}{Acknowledgments} 

\newtheorem{Mthm}{Main Theorem}

\theoremstyle{remark}

\newenvironment{proofmthm}{\begin{proof}[Proof of the Main Theorem]}
                   {\end{proof}}

\begin{document}

\bibliographystyle{amsalpha+}
\title[Non-degenerate hypersurfaces]
{Inversion of adjunction for non-degenerate hypersurfaces}
\maketitle

\begin{abstract}
 We prove a precise inversion of adjunction formula for the
 log variety $(\bC^{n+1},X)$, where $X$ is a non-degenerate
 hypersurface. As a corollary, the minimal log discrepancies
 of non-degenerate normal hypersurface singularities are 
 bounded by dimension.
\end{abstract}

\setcounter{section}{-1}


\section{Introduction}


\footnotetext[1]
{2000 Mathematics Subject Classification. Primary: 14B05
Secondary:14M25, 52B20.}

A {\em log variety} $(X,B)$ is a normal variety $X$ endowed with 
an effective $\R$-Weil divisor $B$ such that $K_X+B$ is 
$\R$-Cartier. For each point $P \in X$, the 
{\em minimal log discrepancy} $a(P;X,B)$ is an invariant 
of the singularity of the log variety $(X,B)$ at $P$. 
In connection to the termination of flips in the Minimal Model 
Program, V.V. Shokurov conjectured certain spectral properties 
of minimal log discrepancies, in particular that they are bounded
by the dimension of the variety~\cite{problems}. This is known 
to hold if $\dim X\le 3$ , or if $X$ is a toric variety and $B$ is 
an invariant divisor ~\cite{reid, mrk, kaw,bor,am}.
Our main result adds to this list the case of non-degenerate
hypersurface singularities:

\begin{Mthm}
 Let $0\in X\subset \bC^{n+1}$ be the germ of a 
normal, non-degenerate hypersurface singularity. Then 
\begin{itemize}
  \item[(i)] $a(0;X)=a(0;\bC^{n+1},X)$.
  \item[(ii)] $a(0;X)\le n$, and equality holds 
  if and only if $X$ is nonsingular.
\end{itemize}
\end{Mthm}

The first statement is a {\em precise inversion of adjunction}
for the log variety $(\bC^{n+1},X)$. Inversion of 
adjunction~\cite{3-flips} has been used by V.V. Shokurov in 
his construction of $3$-fold flips, and is conjectured to 
hold for any log variety (see ~\cite{Ko1,SP}). The effective upper 
bound is an immediate corollary.

As for the proof of the theorem, it is enough to consider 
weighted valuations only, since non-degenerate hypersurfaces 
can be resolved by standard toric methods (Lemma~\ref{resln}).
Via the key Lemma~\ref{inters}, suffices to show that the
log discrepancy function attains its minimum inside some proper
cone of the Newton fan associated with the hypersurface $X$ 
(Proposition~\ref{toricmthm}).

\begin{ack}
The author is a Research Fellow of the Japan Society 
for Promotion of Sciences. Partial support by NSF Grant
DMS-9800807 was received at an initial stage.
\end{ack}


\section{Preliminary}


Let $(X,B)$ be a normal variety $X$ endowed
with an effective $\R$-divisor $B$ such that $K+B$ is 
$\R$-Cartier. For a resolution $\mu:\tilde{X}\to X$, there exists a 
unique divisor $\tilde{B}$ on $Y$ such that 
$\mu^*(K_X+B)=K_{\tilde{X}}+\tilde{B}$ and $\tilde{B}=B$ 
on $\tilde{X}\setminus \Exc(\mu)$.
The identity 
$\tilde{B}=\sum_{E \subset \tilde{X}} (1-a(E;X,B))E$ 
associates to each prime divisor $E$ of $\tilde{X}$ a real 
number $a(E;X,B)$, called the {\em log discrepancy} of $E$ with respect
to $(X,B)$. The invariant $a(E;X,B)$ depends only on the valuation
defined by $E$ on the field of rational functions of $X$, with center 
$c_X(E)=\mu(E)$. 

We say that $\mu$ is a {\em log resolution} of
$(X,B)$ if $\Supp(\mu^{-1}(B))\cup\Exc(\mu)$ is a divisor with
simple normal crossings.
The {\it minimal log discrepancy} 
of a log pair $(X,B)$ at a  point $P \in X$ is 
$$
a(P;X,B)=\inf_{c_X(E)=P} a(E;X,B),
$$
where the infimum is taken after all prime divisors on 
resolutions of $X$, having $P$ as a center on $X$ \cite{problems}. 
The log pair $(X,B)$ has only {\em log canonical singularities} 
if $a(E;X,B) \ge 0$ for every valuation $E$ of $X$. Minimal log 
discrepancies are computed as follows: let $(\tilde{X},\tilde{B})$ 
be a log resolution such that $\mu^{-1}(P)$ is a divisor. 
Let $\cup_i E_i$ be the divisor with normal crossings supporting 
$\Exc(\mu)$ and $\tilde{B}$. If $\min_{c_X(E_i)=P } a(E_i;X,B)<0$
then $a(P;X,B)=-\infty$. Otherwise, 
$
a(P;X,B)= \min_{c_X(E_i)=P } a(E_i;B) \in \R_{\ge 0}
$
and $(X,B)$ has only log canonical 
singularities in a neighborhood of $P$.

\par
To any hypersurface $X:(f=0)\subset {\mathbb C}^{n+1}$ 
one can associate a fan, which is a subdivision of the
standard cone. We recall below this construction, and we also
fix the notation.
Fixing coordinates, we indentify $\bC^{n+1}$ 
with the toric variety $T_N(\sigma)$, where $\sigma$ is the 
standard cone, with basis $\{e_0,\ldots, e_n\}$. 
The {\em Newton polyhedron of} $f$, denoted 
$\Gamma_+$, is the convex hull of $\cup_{m\in \Supp(f)}
(m+\sigma^\vee) \subset M_{\mathbb R}$. 
The {\em Newton diagram} $\Gamma$ is the union of compact 
faces of $\Gamma_+$. The {\em supporting function} 
$l_\Gamma:\sigma \to [0,\infty)$
is defined as
$
l_\Gamma(a)=\min_{m\in \Gamma_+} (a,m).
$
The {\em trace} of a covector $a\in \sigma$ is 
$\tr(a)=\{m\in \Gamma_+; (a,m)=l_\Gamma(a)\}$. The faces 
of $\Gamma_+$ are the traces of covectors $a\in \sigma$. 
Compact faces are the traces of covectors $a\in \relint(\sigma)$, 
while non-compact faces can be written as
$\gamma+ \sigma_I^\vee$, where $\gamma$ is a compact 
face and $\sigma_I^\vee=
\{m\in \sigma^\vee; m_i=0\  \forall i\notin I\}$.

We say that two covectors $a$ and $a'$ are equivalent if they 
have the same trace. The closures of equivalence classes are 
closed cones forming a {\em fan} 
$\Sigma_f=\Sigma_{\Gamma(f)}$, which is a subdivision 
of the standard fan $\{\sigma\}$. Each cone can be represented as
$$
\sigma_{\gamma,I}:=\{a\in \sigma_I; (a,m)=l_\Gamma(a) \ \forall m\in 
\gamma\},\  \sigma_{\gamma,I} \cap \relint(\sigma_I)\ne \emptyset
$$
where $\gamma$ is a compact face of $\Gamma$ and $I$ is a subset 
of $\{1,\ldots,n\}$. We drop $I$ from notation if $I=\emptyset$.
The cones containing $\sigma_{\gamma,I}$ are $\{\sigma_{\tau,J};
\tau\prec \gamma, \ J\subset I\}$. We say that a cone of 
$\Sigma_f$ is {\em proper} if it is not maximal dimensional.
Note that the supporting function $l_\Gamma$ is linear on 
each cone: $l_\Gamma(a)=(a,m)$, if $a\in \sigma_\gamma$ 
and $m$ is a point of $\gamma$.

\begin{defn} The power series $f= \sum_m c_m x^m$ is 
{\em non-degenerate} (with respect to its Newton 
polyhedron) if the hypersurfaces 
$$\{\sum_{m\in \gamma} c_m x^m=0\}
\subset ({\mathbb C}\setminus 0)^{n+1}$$ are non-singular 
for every compact face $\gamma$ of $\Gamma$.
\end{defn}

\begin{lem}\cite[8.9]{AGV} \label{resln}
Let $X \subset {\mathbb C}^{n+1}$ 
be a hypersurface given by a non-degenerate series $f$, and let 
$\Sigma$ be a simple subdivision of $\Sigma_{\Gamma(f)}$ 
containing the primitive vector $\1=(1,\cdots,1)$ in its skeleton.
Then the induced birational morphism $\mu:T_N(\Sigma) \to 
{\mathbb C}^n$ is a log resolution of $({\mathbb C}^n,X)$ over 
a neighbourhood of $0$, and $\mu^{-1}(0)$ is a divisor.
\end{lem}


\section{The log discrepancy function}


For $\delta\in \Shed(\sigma^\vee)$,
the {\em log discrepancy function} is defined as
$$
\varphi=\varphi_\delta: N \cap \sigma \to \R, \ 
\varphi (a)= (a,\delta)-l_{\Gamma}(a)
$$
The function $\varphi$ measures the singularities
of the pair
$$(\bC^{n+1},X+\sum_{i=1}^n (1-\delta_i)H_i)$$ 
where $H_i$
are the coordinate hyperplane sections: if 
$a\in \sigma\cap N$ is a primitive covector, and $E_a$
is the exceptional divisor of the $a$-weighted blow up, then
$a(E_a; \bC^{n+1},X+\sum_{i=1}^n (1-\delta_i)H_i)=
\varphi (a)$. 

We assume until the end of this section that $\varphi$ is 
{\em non-negative}, which is equivalent to 
$\delta \in \Gamma_+ \cap \Shed(\sigma^\vee)$. In 
particular, $\varphi$ has a minimum.

\begin{prop}\label{toricmthm} 
The restriction $\varphi: N \cap \relint(\sigma) \to \R$ 
attains its minimum inside some proper cone of $\Sigma_f$. 
\end{prop}

We first remark that $\varphi$ has the following properties:
\begin{itemize}
   \item[-] $\varphi(ca)=c\cdot \varphi(a)$ for $c>0$.
   \item[-] $\varphi(a+a')\le \varphi(a)+\varphi(a')$, and 
                     equality holds iff $\tr(a)\cap \tr(a')\ne \emptyset$.
    \item[-] The zero locus $\cZ(\varphi):=
     \{a\in \sigma; \varphi(a)=0\}$ is a cone 
	 $\sigma_{\gamma,I}$ of $\Sigma_\Gamma$. 
	 Indeed, one can write (not uniquely)
$
\delta=\sum_i \lambda_i m^i +r,
$
where $\{m^i\}$ are the vertices of a compact face 
$\gamma$ of $\Gamma$, $\lambda_i>0 \ \forall i, 
\sum_i \lambda_i=1$, and $r \in \sigma^\vee$. Then 
$
Z(\varphi)=\sigma_\gamma \cap r^\perp.
$ 
Note that $\delta\in \relint(\gamma+\sigma_I^\vee)$. 
 \item[-] The function $\varphi$ can attain its minimum 
 value only on cones containing $\cZ(\varphi)$. Indeed, 
 there exists $e\in \sigma_{\gamma,I}=\cZ(\varphi)$ 
 such that $\tr_\Gamma(e)=\gamma$. If 
 $a\in \sigma\setminus \cup_{m\in \gamma} \sigma_m$, 
 i.e. $\tr_\Gamma(a)\cap \gamma=\emptyset$, then
 $\varphi(a+e)<\varphi(a)+\varphi(e)=\varphi(a)$. 
Therefore $\varphi$ cannot attain the minimum value at $a$.
\end{itemize}

\begin{lem}\label{keylemma} 
Let $a\in \sigma\cap N$ such that $\tr_\Gamma(a) \cap 
\tr_\Gamma(e_j)=\emptyset$. 
Then $\varphi(a+e_j)\le \varphi(a)$, and equality holds if 
and only if $\delta_j=1$ and one of the following holds:
 \begin{itemize}
   \item[a)] There exists a vertex $m$ of $\Gamma$ 
   such that $a,a+e_j\in \sigma_m$, and $m_j =1$. Note 
   that $\min \varphi|_{\sigma_m}$ is attained on 
    $\cup_{m'_j=0} \sigma_{mm'}$ in this case.      
   \item[b)] $m'_j=0$ and $(a,m')=(a,m)+1$ for every
   $m\in \tr_\Gamma(a)$ and $m'\in \tr_\Gamma(a+e_j)$. 
   In particular, $e_j \in \sigma_{m'}$.
  \end{itemize}
\end{lem}

\begin{proof} 
Let $m\in \tr_\Gamma(a), m'\in \tr_\Gamma(a+e_j)$. Then
$$
\varphi(a+e_j)-\varphi(a)=-(a,m'-m)+\delta_j-m'_j.
$$
 \begin{itemize}
 \item[a)] If there exists 
   $m\in \tr_\Gamma(a) \cap \tr_\Gamma(a+e_j)$, then
   $\varphi(a+e_j)-\varphi(a)=\delta_j-m_j$. Since 
   $e_j\notin \sigma_m$, $m_j\ge 1$, and there exists 
   $m'$ with $m'_j<m_j$. Thus $\varphi(a+e_j) \le 
   \varphi(a)$, and equality holds if and only if
   $\delta_j=m_j=1$.
  \item[b)] 
   Assume $\tr_\Gamma(a) \cap \tr_\Gamma(a+e_j)=\emptyset$. 
   Then $\varphi(a+e_j)-\varphi(a)\le -1+\delta_j-m'_j\le 0$, 
   and equality holds if and only if $\delta_j=1$, $m'_j=0$, and
   $(a,m'-m)=1$ for all $m\in \tr_\Gamma(a), m'\in 
   \tr_\Gamma(a+e_j)$.
  \end{itemize}
\end{proof}

\begin{proof}[Proof of Proposition \ref{toricmthm}]
Assume first that $I=\emptyset$. That is 
$\delta\in \Gamma$, and $\cZ(\varphi)=\sigma_\gamma$ 
intersects $\relint(\sigma)$. If $\dim \gamma>0$ then 
the minimum is attained only in the proper cone 
$\sigma_\gamma$.
If $\dim \gamma =0$, $\varphi$ is identically zero on 
the maximal cone $\sigma_\gamma$.
\par
Assume $I\ne \emptyset$. We may assume 
$\delta_i<1$ for all $i\in I$. This is sufficient, since 
the desired property of $\varphi=\varphi_\delta$ is 
closed with respect to $\delta$ belonging to the convex 
polyhedron $\Gamma_+\cap \Shed(\sigma^\vee)$. Assume 
by contradiction that minimum is not attained on proper 
cones. Let $m$ be a vertex of $\gamma$ such that $\varphi$ 
attains its minimum inside the maximal cone $\sigma_m$. 
Since $\delta_i<1$ for all $i\in I$, we obtain 
$e_i\in \sigma_m$ for all $i\in I$ by Lemma~\ref{keylemma}. 
In fact, we obtain $m_i=0$ for all $i\in I$ since $\varphi$ is 
non-negative and $\delta_i<1$. By Lemma~\ref{keylemma} 
again, $\varphi$ attains minimum inside $\sigma_{m'}$ for 
some $m'\ne m$.

Let $m\in \gamma$ with $m_i=0 \ \forall i\in I$. Then the
minimum of $\varphi$ on $\sigma_m\cap \relint(\sigma)$ 
is attained on $\sigma_{\gamma,I}+\sum_{i\in I}e_i$. 
Indeed, for $a\in \sigma_m$ we have
$$
\varphi(a)=\sum_{i\in I}\delta_i a_i+ \varphi(\bar{a})\ge 
\sum_{i\in I}\delta_i, 
$$
where $\bar{a}_i=a_i$ for $i\notin I$ and $\bar{a}_i=0$ 
for $i\in I$. If $m,m'\in \gamma$ with $m_i=m'_i=0 \ 
\forall i\in I$ then $\sigma_{\gamma,I}+\sum_{i\in I}
e_i \subset \sigma_{mm'}$ and we are done.
\end{proof}


\section{Inversion of Adjunction}


\begin{lem}\label{inters} 
Let $X:(f(x)=0)\subset {\mathbb C}^{n+1}$ be a 
hypersurface which does not contain any of the coordinate
hyperplanes $(x_i=0)$. Let $\Sigma$ be a simple subdivision 
of $\Sigma_{\Gamma(f)}$, and let $\mu:T_N(\Sigma) \to 
{\mathbb C}^{n+1}$ be the associated resolution.
Let $E_a$ be the $\mu$-exceptional divisor corresponding to 
$a\in \relint(\sigma)\cap N$. If $E_a$ does not intersect the 
proper transform of $X$, then $a$ belongs to a unique cone of 
$\Sigma_{\Gamma(f)}$.
\end{lem}

\begin{proof}
$T_N(\Sigma)$ is covered by the open sets 
$U_\tau \simeq {\bC}^{n+1}$ corresponding to the maximal 
dimensional cones of $\Sigma$. Let $\tau \in \Sigma$ be a 
maximal cone such that $E_a\cap U_\tau \ne \emptyset$. 
The restriction $\mu:U_\tau \to \bC^{n+1}$ can be identified 
with
$$
\mu_\tau : \bC^{n+1} \to  \bC^{n+1}, \ \ 
x_i=y_0^{a^0_i} \cdots y_n^{a^n_i}
$$
where $(a^0,\ldots, a^n)$ is the ordered skeleton of $\tau$.
We we may assume $a=a^0$. Denote by $X'$ the proper 
transform of $X$ in $T_N(\Sigma)$. 
We have
$$
{\mu_\tau}^*(f)=
\sum_{m\in \Supp(f)} c_m \prod_{j=0}^n y_j^{(a^j,m)}=
(\prod_{j=0}^n y_j^{l(a^j)}) f_\tau(y_0,\ldots,y_n),
$$
where $f_\tau=0$ is the equation of $X'\cap U_\tau$. 
The divisor $E_a$ has equation $y_0=0$ in $U_\tau$, thus 
$E_a\cap X'\cap U_\tau=\emptyset$ iff there is a 
non-zero constant $C$ such that
$f_\tau \equiv C \mod y_0$.
Equivalently, $\min_{m\in \Supp(f)} (a^0,m)$ is attained in 
exactly one point. Since $a\in \relint(\sigma)$, this implies 
that the trace of $a$ is a vertex of $\Gamma$, i.e. $a$ belongs 
to the interior of some maximal cone of $\Sigma_\Gamma$.
\end{proof}

\begin{proofmthm} (i):
We may assume that $X$ does 
not contain any of the coordinate hyperplanes. 
By Lemma~\ref{resln}, we obtain 
$$
a(0;X)\ge a(0;\bC^{n+1},X)=\inf 
\varphi|_{\relint(\sigma)\cap N}
$$
where $\varphi$ corresponds to $\delta=\1$. 
Assume that $a(0;\bC^{n+1},X)=-\infty$. Then $\varphi$
takes negative values. In particular, there exists a
primitive vector $a\in \relint(\sigma)\cap N$ contained
in a proper cone of $\Sigma_f$ such that $\varphi(a)<0$.
Let $\Sigma$ be a simple subdivision of the fan 
$\Sigma_\Gamma$, such that $a$ and $\1$ belong to 
its skeleton, and let $\mu:T_N(\Sigma)\to {\mathbb C}^n$ 
be the induced log resolution. By Lemma~\ref{inters},
the exceptional divisor $E_a$ intersects the proper 
transform $X'$ of $X$, and Lemma~\ref{resln} 
implies that $a(0;X)\le a(E_a\cap X';X)=\varphi(a)<0$. 
Therefore $a(0;X)=a(0;\bC^{n+1},X)$.

Assume that $a(0;\bC^{n+1},X)\ge 0$. In particular, $\varphi$
is non-negative. By Proposition~\ref{toricmthm}, there exists 
a primitive vector $a\in \relint(\sigma)\cap N$ contained in
a proper cone of $\Sigma_f$ such that 
$\varphi(a)=\min \varphi|_{\relint(\sigma)\cap N}$.
The same argument as above implies that 
$a(0;X)\le \varphi(a)$, hence $a(0;X)=a(0;\bC^{n+1},X)$.

(ii): If $E$ is the exceptional divisor of the blow up of
$\bC^{n+1}$ at $0$, then $a(E;\bC^{n+1},X)=n-(\mult_0(f)-1)$.
By $(i)$, $a(0;X)=a(0;\bC^{n+1},X) \le n$. It is clear that equality
holds if and only if $X$ is non-singular at $0$.
\end{proofmthm}


\end{document}